\documentclass[a4paper, 11pt, twoside]{article}
\usepackage{amsmath}
\usepackage{amscd}
\usepackage{amssymb}
\usepackage{amsthm}
\usepackage{graphics}
\usepackage{indentfirst}
\usepackage{latexsym}
\usepackage{amsfonts}
\usepackage{multicol}
\usepackage{stmaryrd}
\usepackage{array}
\usepackage{youngtab}
\usepackage{pxfonts}

\newtheoremstyle{mattthm}{}{}{\itshape}{}{\bfseries}{.}{ }{}

\theoremstyle{mattthm}
\newtheorem{lemma}{Lemma}[section]
\newtheorem{propn}[lemma]{Proposition}
\newtheorem{thm}[lemma]{Theorem}

\newtheorem{conj}[lemma]{Conjecture}

\newtheoremstyle{mattdef}{}{}{}{}{\bfseries}{.}{ }{}

\theoremstyle{mattdef}
\newtheorem*{rmk}{Remark}

\newtheorem*{eg}{Example}
\newtheorem*{egs}{Examples}

\begin{document}

\newenvironment{pf}{\noindent\textbf{Proof.}}{\hfill \qedsymbol\newline}
\newenvironment{pfof}[1]{\vspace{\topsep}\noindent\textbf{Proof of {#1}.}}{\hfill \qedsymbol\newline}
\newenvironment{pfenum}{\noindent\textbf{Proof.}\indent\begin{enumerate}\vspace{-\topsep}}{\end{enumerate}\vspace{-\topsep}\hfill \qedsymbol\newline}
\newenvironment{pfnb}{\noindent\textbf{Proof.}}{\newline}

\newcommand\la\lambda
\newcommand{\sss}{\mathfrak{S}_}
\newcommand\sgn{\operatorname{sgn}}
\newcommand{\thmcite}[2]{\textup{\textbf{\cite[#2]{#1}}}\ }
\newcommand{\bbf}{\mathbb{F}}
\newcommand{\bbn}{\mathbb{N}}
\newcommand{\bbq}{\mathbb{Q}}
\newcommand\calu{\mathcal{U}}
\newcommand\calf{\mathcal{F}}
\newcommand\nchar{\operatorname{char}}
\newcommand\lad{\operatorname{lad}}
\newcommand{\gs}{\geqslant}
\newcommand{\ls}{\leqslant}
\newcommand{\sect}[1]{\section{#1}}
\newcommand{\clam}{\begin{description}\item[\hspace{\leftmargin}Claim.]}
\newcommand{\prof}{\item[\hspace{\leftmargin}Proof.]}
\newcommand{\malc}{\end{description}}

\title{Regularisation and the Mullineux map}
\author{Matthew Fayers\footnote{This research was undertaken while the author was visiting Massachusetts Institute.of Technology as a Postdoctoral Fellow, with the support of a Research Fellowship from the Royal Commission for the Exhibition of 1851; the author is very grateful to M.I.T. for its hospitality, and to the 1851 Commission for its generous support.}\\\normalsize  Queen Mary, University of London, Mile End Road, London E1 4NS, U.K.\\\texttt{\normalsize  m.fayers@qmul.ac.uk}}
\date{}
\maketitle
\begin{center}
2000 Mathematics subject classification: 05E10, 20C30
\end{center}
\markboth{Matthew Fayers}{Regularisation and the Mullineux map}
\pagestyle{myheadings}

\begin{abstract}
We classify the pairs of conjugate partitions whose regularisations are images of each other under the Mullineux map.  This classification proves a conjecture of Lyle, answering a question of Bessenrodt, Olsson and Xu.
\end{abstract}

\sect{Introduction}

Suppose $n\gs0$ and $\bbf$ is a field of characteristic $p$; we adopt the convention that the characteristic of a field is the order of its prime subfield.  It is well known that the representation theory of the symmetric group $\sss n$ is closely related to the combinatorics of partitions.  In particular, for each partition $\la$ of $n$, there is an important $\bbf\sss n$-module $S^\la$ called the \emph{Specht module}.  If $p=\infty$, then the Specht modules are irreducible and afford all irreducible representations of $\bbf\sss n$.  If $p$ is a prime, then for each $p$-regular partition $\la$ the Specht module $S^\la$ has an irreducible cosocle $D^\la$, and the modules $D^\la$ afford all irreducible representations of $\bbf\sss n$ as $\la$ ranges over the set of $p$-regular partitions of $n$.

Given this set-up, it is natural to express representation-theoretic statements in terms of the combinatorics of partitions.  An example of this which is of central interest in this paper is the \emph{Mullineux map}.  Let $\sgn$ denote the one-dimensional sign representation of $\bbf\sss n$.  Then there is an involutory functor $-\otimes\sgn$ from the category of $\bbf\sss n$-modules to itself.  This functor sends simple modules to simple modules, and therefore for each $p$-regular partition $\la$ there is some $p$-regular partition $M(\la)$ such that $D^\la\otimes\sgn\cong D^{M(\la)}$.  The map $M$ thus defined is now called the Mullineux map, since it coincides with a map defined combinatorially by Mullineux \cite{mull}; this was proved by Ford and Kleshchev \cite{fk}, using an alternative combinatorial description of $M$ due to Kleshchev \cite{k3}.

Another important aspect of the combinatorics of partitions from the point of view of representation theory is \emph{$p$-regularisation}.  This combinatorial procedure was defined by James in order to describe, for each partition $\la$, a $p$-regular partition (which is denoted $G\la$ in this paper) such that the simple module $D^{G\la}$ occurs exactly once as a composition factor of $S^\la$.  In this paper we study the relationship between the Mullineux map and regularisation.  Our motivation is the observation that if $p=2$ or $p$ is large relative to the size of $\la$, then $MG\la=GT\la$, where $T\la$ denotes the conjugate partition to $\la$.  However, this is not true for arbitrary $p$, and it natural to ask for which pairs $(p,\la)$ we have $MG\la=GT\la$.  The purpose of this paper is to answer this question, which was first posed by Bessenrodt, Olsson and Xu; the answer confirms a conjecture of Lyle.

If we replace the group algebra $\bbf\sss n$ with the Iwahori--Hecke algebra of the symmetric group at a primitive $e$th root of unity in $\bbf$ (for some $e\gs2$), then all of the above background holds true, with the prime $p$ replaced by the integer $e$ (and with an appropriate analogue of the sign representation).  Therefore, in this paper, we work with an arbitrary integer $e\gs2$ rather than a prime $p$.

In the remainder of this section we give all the definitions we shall need concerning partitions, and state our main result.  Section \ref{if} is devoted to proving one half of the conjecture, and Section \ref{onlyif} to the other half.  While the first half of the proof consists of elementary combinatorics, the latter half of the proof is algebraic, being an easy consequence of two theorems about $v$-decomposition numbers in the Fock space.  We introduce the background material for this as we need it.

\subsection{Partitions}\label{backsec}

A \emph{partition} is a sequence $\la=(\la_1,\la_2,\dots)$ of non-negative integers such that $\la_1\gs\la_2\gs\dots$ and the sum $|\la|=\la_1+\la_2+\dots$ is finite.  We say that $\la$ is a partition of $|\la|$.  When writing partitions, we usually group together equal parts and omit zeroes.  We write $\varnothing$ for the unique partition of $0$.

$\la$ is often identified with its \emph{Young diagram}, which is the subset
\[[\la]=\left\{(i,j)\mid j\ls\la_i\right\}\]
of $\bbn^2$.  We refer to elements of $\bbn^2$ as \emph{nodes}, and to elements of $[\la]$ as nodes of $\la$.  We draw the Young diagram as an array of boxes using the English convention, so that $i$ increases down the page and $j$ increases from left to right.

If $e\gs2$ is an integer, we say that $\la$ is \emph{$e$-regular} if there is no $i\gs1$ such that $\la_i=\la_{i+e-1}>0$, and otherwise we say that $\la$ is \emph{$e$-singular}.  We say that $\la$ is \emph{$e$-restricted} if $\la_i-\la_{i+1}<e$ for all $i\gs1$.


\subsection{Operators on partitions}

Here we introduce a variety of operators on partitions.  These include regularisation and the Mullineux map, as well as other more familiar operators which will be useful.

\subsubsection{Conjugation}

Suppose $\la$ is a partition.  The \emph{conjugate partition} to $\la$ is the partition $T\la$ obtained by reflecting the Young diagram along the main diagonal.  That is,
\[(T\la)_i = \left|\left\{j\gs1\ \left|\ \la_j\gs i\right.\right\}\right|.\]
We remark that $T\la$ is conventionally denoted $\la'$; we choose our notation in this paper so that all operators on partitions are denoted with capital letters written on the left.  The letter $T$ is taken from \cite{box}, and stands for `transpose'.

In this paper we write $l(\la)$ for $(T\la)_1$, i.e.\ the number of non-zero parts of $\la$. 

\subsubsection{Row and column removal}

Suppose $\la$ is a partition.  Let $R\la$ denote the partition obtained by removing the first row of the Young diagram; that is, $(R\la)_i=\la_{i+1}$ for $i\gs1$.  Similarly, let $C\la$ denote the partition obtained by removing the first column from the Young diagram of $\la$, i.e.\ $(C\la)_i = \max\{\la_i-1,0\}$ for $i\gs1$.

In this paper we shall use without comment the obvious relation $TR=CT$.

\subsubsection{Regularisation}

Now we introduce one of the most important concepts of this paper.  Suppose $\la$ is a partition and $e\gs2$.  The \emph{$e$-regularisation} of $\la$ is an $e$-regular partition associated to $\la$ in a natural way.  The notion of regularisation was introduced by James \cite{j1} in the case where $e$ is a prime, where it plays a r\^ ole in the computation of the $e$-modular decomposition matrices of the symmetric groups.

For $l\gs1$, we define the $l$th \emph{ladder} in $\bbn^2$ to be the set of nodes $(i,j)$ such that $i+(e-1)(j-1)=l$.  The regularisation of $\la$ is defined by moving all the nodes of $\la$ in each ladder as high as they will go within that ladder.  It is a straightforward exercise to show that this procedure gives the Young diagram of a partition, and the $e$-regularisation of $\la$ is defined to be this partition.

\begin{eg}
Suppose $e=3$, and $\la=(4,3^3,1^5)$.  Then the $e$-regularisation of $\la$ is $(5,4,3^2,2,1)$, as we can see from the following Young diagrams, in which we label each node with the number of the ladder in which it lies.
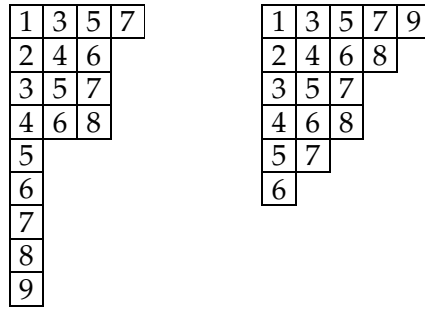
\[\young(1357,246,357,468,5,6,7,8,9)\qquad\qquad\raisebox{37.8pt}{\young(13579,2468,357,468,57,6)}\]
\end{eg}

We write $G\la$ for the $e$-regularisation of $\la$.  Clearly $G\la$ is $e$-regular, and equals $\la$ if $\la$ is $e$-regular.  We record here three results we shall need later; the proofs of the first two are easy exercises.

\begin{lemma}\label{rggr}
Suppose $\la$ is a partition.  If $(G\la)_1=\la_1$, then $RG\la=GR\la$.
\end{lemma}

\begin{lemma}\label{cggc}
Suppose $\la$ and $\mu$ are partitions.  If $l(\la)=l(\mu)$ and $GC\la=C\mu$, then $G\la=G\mu$.
\end{lemma}

\begin{lemma}\label{simplereg}
Suppose $\zeta$ is an $e$-regular partition, and $x\gs l(\zeta)+e-1$.  Let $\xi$ be the partition obtained by adding a column of length $x$ to $\zeta$, and let $\eta$ be the partition obtained by adding a column of length $x-e+1$ to $C\zeta$.  Then $G\eta=CG\xi$.
\end{lemma}

\begin{pf}
For any $n\gs1$ and any partition $\la$, let $\lad_n(\la)$ denote the number of nodes of $\la$ in ladder $n$.  Since $G\eta$ and $CG\xi$ are both $e$-regular, it suffices to show that $\lad_n(G\eta)=\lad_n(CG\xi)$ for all $n$.

$\eta$ is obtained from $\zeta$ by adding the nodes $(l(\zeta)+1,1),\dots,(x-e+1,1)$, so we have
\[\lad_n(G\eta) = \lad_n(\eta) = \begin{cases}
\lad_n(\zeta)+1 & (l(\zeta)<n< x+e)\\
\lad_n(\zeta) & (\text{otherwise}).
\end{cases}\]
It is also easy to compute
\[\lad_n(\xi) = \begin{cases}
1 & (1\ls n<e)\\
\lad_{n-e+1}(\zeta)+1 & (e\ls n\ls x)\\
\lad_{n-e+1}(\zeta) & (x<n).
\end{cases}\]
\clam
$l(G\xi) = l(\zeta)+e-1$.
\prof
Since $\zeta$ is $e$-regular and $(l(\zeta),1)\in[\zeta]$, every node of ladder $l(\zeta)$ is a node of $\zeta$.  Hence every node of ladder $l(\zeta)+e-1$ is a node of $\xi$; so when $\xi$ is regularised, none of these nodes moves, and we have $(l(\zeta)+e-1,1)\in[G\xi]$, i.e.\ $l(G\xi)\gs l(\zeta)+e-1$.

On the other hand, the node $(l(\zeta)+1,2)$ does not lie in $[\xi]$, so the node $(l(\zeta)+e,1)$ cannot lie in $[G\xi]$, i.e.\ $l(G\xi)<l(\zeta)+e$.
\malc

From the claim we deduce that
\[\lad_n(CG\xi) = \begin{cases}
\lad_{n+e-1}(\xi)-1 & (n\ls l(\zeta))\\
\lad_{n+e-1}(\xi) & (n>l(\zeta)),
\end{cases}\]
and combining this with the statements above gives the result.
\end{pf}

\subsubsection{The Mullineux map}

Now we introduce the Mullineux map, which is the most important concept of this paper.  
We shall give two different recursive definitions of the Mullineux map: the original definition due to Mullineux \cite{mull}, and an alternative version due to Xu \cite{xu}.

Suppose $\la$ is a partition, and define the \emph{rim} of $\la$ to be the subset of $[\la]$ consisting of all nodes $(i,j)$ such that $(i+1,j+1)\notin\la$.  Now fix $e\gs2$, and suppose that $\la$ is $e$-regular.  Define the \emph{$e$-rim} of $\la$ to be the subset $\left\{(i_1,j_1),\dots,(i_r,j_r)\right\}$ of the rim of $\la$ obtained by the following procedure.

\begin{itemize}
\item
If $\la=\varnothing$, then set $r=0$, so that the $e$-rim of $\la$ is empty.  Otherwise, let $(i_1,j_1)$ be the top-rightmost node of the rim, i.e.\ the node $(1,\la_1)$.
\item
For $k>1$ with $e\nmid k-1$, let $(i_k,j_k)$ be the next node along the rim from $(i_{k-1},j_{k-1})$, i.e.\ the node $(i_{k-1}+1,j_{k-1})$ if $\la_{i_{k-1}}=\la_{i_{k-1}+1}$, or the node $(i_{k-1},j_{k-1}-1)$ otherwise.
\item
For $k>1$ with $e\mid k-1$, define $(i_k,j_k)$ to be the node $(i_{k-1}+1,\la_{i_{k-1}+1})$.
\item
Continue until a node $(i_k,j_k)$ is reached in the bottom row of $[\la]$ (i.e.\ with $i_k=l(\la)$), and either $j_k=1$ or $e\mid k$.  Set $r=k$, and stop.
\end{itemize}

Less formally, we construct the $e$-rim of $\la$ by working along the rim from top right to bottom left, and moving down one row every time the number of nodes we've seen is divisible by $e$.

The integer $r$ defined in this way is called the \emph{$e$-rim length} of $\la$.  We define $I\la$ to be the partition obtained by removing the $e$-rim of $\la$ from $[\la]$.

\begin{egs}\indent\vspace{-\topsep}
\begin{enumerate}
\item
Suppose $e=3$, and $\la=(10,6^2,4,2)$.  Then the $e$-rim of $\la$ consists of the marked nodes in the following diagram, and we see that $r=11$ and $I\la=(7,5,4,1)$.
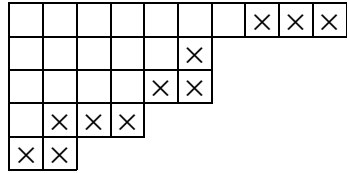
\[\young(\ \ \ \ \ \ \ \times\times\times,\ \ \ \ \ \times,\ \ \ \ \times\times,\ \times\times\times,\times\times)\]
\item
Suppose $e=2$, and $\la$ is any $2$-regular partition.  The $2$-rim of $\la$ consists of the last two nodes in each row of $[\la]$ (or the last node, if there is only one).  Hence when $e=2$ the operator $I$ is the same as $C^2$.
\end{enumerate}
\end{egs}

Now we can define the Mullineux map recursively.  Suppose $\la$ is an $e$-regular partition.  If $\la=\varnothing$, then set $M\la=\varnothing$.  Otherwise, compute the partition $I\la$ as above.  Then $|I\la|<|\la|$, and $I\la$ is $e$-regular, so we may assume that $MI\la$ is defined.  Let $r$ be the $e$-rim length of $\la$, and define
\[m = \begin{cases}
r-l(\la) & (e\mid r)\\
r-l(\la)+1 & (e\nmid r).
\end{cases}\]
It turns out that there is a unique $e$-regular partition $\mu$ which has $e$-rim length $r$ and $l(\mu)=m$, and which satisfies $I\mu=MI\la$.  We set $M\la=\mu$.

\begin{egs}\indent\vspace{-\topsep}
\begin{enumerate}
\item
Suppose $e=3$, $\la=(3^2,2^2,1)$ and $\mu=(6,4,1)$.  Then we have $I\la=(2,1^2)$ and $I\mu=(3,1)$, as we see from the following diagrams.
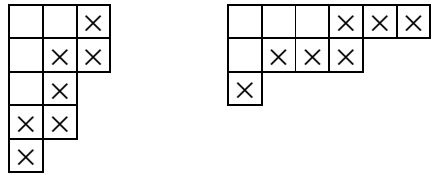
\[\young(\ \ \times,\ \times\times,\ \times,\times\times,\times)\qquad\qquad\raisebox{25.25pt}{\young(\ \ \ \times\times\times,\ \times\times\times,\times)}\]
Computing $e$-rims again, we find that $I^2\la=I^2\mu=\varnothing$.  Now comparing the numbers of non-zero parts of these partitions with their $e$-rim lengths we find that $MI\la=I\mu$, and hence that $M\la=\mu$.
\item
Suppose $e=2$, and $\la$ is a $2$-regular partition.  From above, we see that the $2$-rim length of $\la$ is $2l(\la)$, if $\la_{l(\la)}\gs2$, or $2l(\la)-1$ if $\la_{l(\la)}=1$.  Either way, we get $m=l(\la)$, and this implies inductively that in the case $e=2$ the Mullineux map is the identity.
\item
Suppose $e$ is large relative to $\la$; in particular, suppose $e$ is greater than the number of nodes in the rim of $\la$.  Then the $e$-rim of $\la$ coincides with the rim, so that the $e$-rim length is $\la_1+l(\la)-1$.  Hence $m=\la_1$, and from this it is easy to prove by induction that $M\la=T\la$.
\end{enumerate}
\end{egs}

Now we give Xu's alternative definition of the Mullineux map.  Suppose $\la$ is a partition with $e$-rim length $r$, and define
\[l' = \begin{cases}
l(\la) & (e\mid r)\\
l(\la)-1 & (e\nmid r).
\end{cases}\]
Define $J\la$ to be the partition obtained by removing the $e$-rim from $\la$, and then adding a column of length $l'$.  Another way to think of this is to define the \emph{truncated $e$-rim} of $\la$ to be the set of nodes $(i,j)$ in the $e$-rim of $\la$ such that $(i,j-1)$ also lies in the $e$-rim, together with the node $(l(\la),1)$ if $e\nmid r$, and to define $J\la$ to be the partition obtained by removing the truncated $e$-rim.

\begin{eg}
Returning to an earlier example, take $e=3$ and $\la=(10,6^2,4,2)$.  Then the truncated $e$-rim of $\la$ consists of the marked nodes in the following diagram, and we see that $I\la=(8,6,5,2)$.
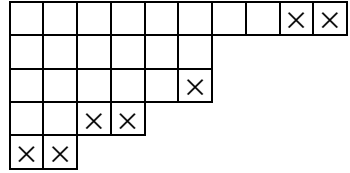
\[\young(\ \ \ \ \ \ \ \ \times\times,\ \ \ \ \ \ ,\ \ \ \ \ \times,\ \ \times\times,\times\times)\]
\end{eg}

If $\la$ is $e$-regular, then it is a simple exercise to show that $J\la$ is $e$-regular and $|J\la|<|\la|$.  So we assume that $MJ\la$ is defined recursively, and we define $M\la$ to be the partition obtained by adding a column of length $|\la|-|J\la|$ to $MJ\la$.  Xu \cite[Theorem 1]{xu} shows that this map coincides with Mullineux's map $M$.  In other words, we have the following.

\begin{propn}\label{xuv}
Suppose $\la$ and $\mu$ are $e$-regular partitions, with $|\la|=|\mu|$.  Then $M\la=\mu$ if and only if $MJ\la=C\mu$.
\end{propn}

\subsection{Hooks}\label{hooks}

Now we set up some basic notation concerning hooks in Young diagrams.  Suppose $\la$ is a partition, and $(i,j)$ is a node of $\la$.  The \emph{$(i,j)$-hook} of $\la$ is defined to be the set $H_{ij}(\la)$ of nodes in $[\la]$ directly to the right of or directly below $(i,j)$, including the node $(i,j)$ itself.  The \emph{arm length} $a_{ij}(\la)$ is the number of nodes directly to the right of $(i,j)$, i.e.\ $\la_i-j$, and the \emph{leg length} $l_{ij}(\la)$ is the number of nodes directly below $(i,j)$, i.e.\ $(T\la)_j-i$.  The \emph{$(i,j)$-hook length} $h_{ij}(\la)$ is the total number of nodes in $H_{ij}(\la)$, i.e.\ $a_{ij}(\la)+l_{ij}(\la)+1$.

Now fix $e\gs2$.  The \emph{$e$-weight} of $\la$ is defined to be the number of nodes $(i,j)$ of $\la$ such that $e\mid h_{ij}(\la)$.  If $(i,j)\in[\la]$ with $e\mid h_{ij}(\la)$, we say that $H_{ij}(\la)$ is
\begin{itemize}
\item
\emph{shallow} if $a_{ij}(\la)\gs (e-1)l_{ij}(\la)$, or
\item
\emph{steep} if $l_{ij}(\la)\gs(e-1)a_{ij}(\la)$.
\end{itemize}

\begin{eg}
Suppose $e=3$ and $\la=(5,2,1^4)$.  Then we have $(2,1)\in[\la]$, with $a_{2,1}(\la)=1$, $l_{2,1}(\la)=4$, and hence $h_{2,1}(\la)=6$.  $H_{2,1}(\la)$ is steep if $e=3$, but not if $e=6$.
\end{eg}

\subsection{Lyle's Conjecture}

Suppose $e\gs2$ and $\la$ is an $e$-regular partition.  As noted above, if $e$ is large relative to $|\la|$, then $M\la=T\la$.  Of course, there is no hope that this is true in general, since $T\la$ will not in general be an $e$-regular partition.  But $e$-regularisation provides a natural way to obtain an $e$-regular partition from an arbitrary partition, and it is therefore natural to ask: for which $e$-regular partitions $\la$ do we have $M\la = GT\la$?  When $e$ is large relative to $\la$ we have $G\la=\la$ and (from the example above) $M\la=T\la$, so certainly $M\la=GT\la$ in this case.  We also have $M\la=GT\la$ for all partitions $\la$ when $e=2$: we have seen that for $e=2$ the Mullineux map is the identity, and it is a simple exercise to show that $\la$ and $T\la$ have the same $2$-regularisation for any $\la$.  But it is not generally true that $M\la=GT\la$ for an $e$-regular partition $\la$.  Bessenrodt, Olsson and Xu \cite{box} have given a classification of the partitions for which this does hold, as follows.

\begin{thm}\label{boxthm}\thmcite{box}{Theorem 4.8}
Suppose $\la$ is an $e$-regular partition.  Then $M\la=GT\la$ if and only if for every $(i,j)\in[\la]$ with $e\mid h_{ij}(\la)$, the hook $H_{ij}(\la)$ is shallow.
\end{thm}

\begin{eg}
Suppose $e=4$ and $\la=(14,10,2^2)$.  The Young diagram is as follows; we have marked those nodes $(i,j)$ for which $4\mid h_{ij}(\la)$.
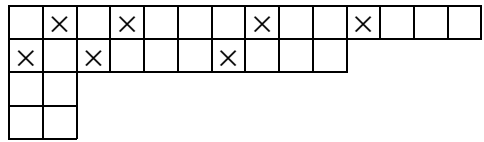
\[\young(\ \times\ \times\ \ \ \times\ \ \times\ \ \ ,\times\ \times\ \ \ \times\ \ \ ,\ \ ,\ \ )\]
We see that all the hooks of length divisible by $4$ are shallow, so $\la$ satisfies the second hypothesis of Theorem \ref{boxthm}.  And it may be verified that $GT\la=M\la=(5^2,4^2,3^2,2^2)$.
\end{eg}

Bessenrodt, Olsson and Xu have also posed the following more general question \cite[p.\ 454]{box}, which is essentially the same problem without the assumption that $\la$ is $e$-regular.
\begin{quote}
For which partitions $\la$ is it true that $MG\la = GT\la$?
\end{quote}
Motivated by the (now solved) problem of the classification of irreducible Specht modules for symmetric groups, Lyle conjectured the following solution in her thesis.

\begin{conj}\label{slconj}\thmcite{slthesis}{Conjecture 5.1.18}
Suppose $\la$ is a partition.  Then $MG\la = GT\la$ if and only if for every $(i,j)\in[\la]$ with $e\mid h_{ij}(\la)$, the hook $H_{ij}(\la)$ is either shallow or steep.
\end{conj}

The purpose of this paper is to prove this conjecture.  It is a simple exercise to show that a partition possessing a steep hook must be $e$-singular; so in the case where $\la$ is $e$-regular, Conjecture \ref{slconj} reduces to Theorem \ref{boxthm}.

Let us define an \emph{L-partition} to be a partition satisfying the second condition of Conjecture \ref{slconj}, i.e.\ a partition for which every $H_{ij}(\la)$ of length divisible by $e$ is either shallow or steep.

\begin{eg}
Suppose $e=4$ and $\la=(11,2^2,1^5)$.  The Young diagram of $\la$ is as follows.
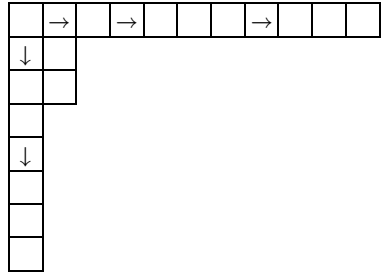
\[\young(\ \shortrightarrow\ \shortrightarrow\ \ \ \shortrightarrow\ \ \ ,\shortdownarrow\ ,\ \ ,\ ,\shortdownarrow,\ ,\ ,\ )\]
The nodes $(i,j)$ with $4\mid h_{ij}(\la)$ are marked; we see that those marked $\young(\shortrightarrow)$ correspond to shallow hooks, and those marked $\young(\shortdownarrow)$ correspond to steep hooks.  So $\la$ is an L-partition when $e=4$.  We have $G\la=(11,3,2^2,1^2)$, $GT\la=(8,4,3^2,2)$, and it can be checked that $MG\la=GT\la$.
\end{eg}


\section{The `if' part of Conjecture \ref{slconj}}\label{if}

In this section we prove the `if' half of Conjecture \ref{slconj}, i.e.\ that $MG\la=GT\la$ whenever $\la$ is an L-partition.  We begin by noting some properties of L-partitions, and making some more definitions.  Note that when $e=2$, every partition is an L-partition; by the above remarks we have $MG\la=GT\la$ for every partition when $e=2$, so Conjecture \ref{slconj} holds when $e=2$.  Therefore, \emph{we assume throughout this section that $e\gs3$}.  The following simple observations will be used without comment.

\begin{lemma}\label{simplesss}
Suppose $\la$ is a partition.  Then $\la$ is an L-partition if and only if $T\la$ is.   If $\la$ is an L-partition, then so are $R\la$ and $C\la$.
\end{lemma}

Now we examine the structure of L-partitions in more detail.  Suppose $\la$ is an L-partition, and let $s(\la)$ be maximal such that $\la_{s(\la)}-\la_{s(\la)+1}\gs e$, setting $s(\la)=0$ if $\la$ is $e$-restricted.  Similarly, set $t(\la)=0$ if $\la$ is $e$-regular, and otherwise let $t(\la)$ be maximal such that $(T\la)_{t(\la)}-(T\la)_{t(\la)+1}\gs e$.  Clearly, we have $s(\la)=t(T\la)$.

\begin{lemma}\label{gse-1}
If $\la$ is an L-partition, then for $1\ls i\ls s(\la)$ we have $\la_i-\la_{i+1}\gs e-1$, while for $1\ls j\ls t(\la)$ we have $(T\la)_j-(T\la)_{j+1}\gs e-1$.
\end{lemma}

\begin{pf}
We prove the first statement.  Suppose this statement is false, and let $i<s(\la)$ be maximal such that $\la_i-\la_{i+1}<e-1$.  Put $j=\la_i-e+2$.  Then we have $(i,j)\in[\la]$, with $a_{ij}(\la)=e-2$ and $l_{ij}(\la)=1$, which (given our assumption that $e\gs 3$) contradicts the assumption that $\la$ is an L-partition.
\end{pf}

\begin{lemma}\label{shallowsteep}
Suppose $\la$ is an L-partition and $(i,j)\in[\la]$ with $e\mid h_{ij}(\la)$.
\begin{enumerate}
\item
If $i>s(\la)$, then $H_{ij}(\la)$ is steep.
\item
If $j>t(\la)$, then $H_{ij}(\la)$ is shallow.
\end{enumerate}
\end{lemma}

\begin{pf}
We prove (1).  Let $a=a_{ij}(\la)$ and $l=l_{ij}(\la)$.  $\la$ is an L-partition, so if $H_{ij}(\la)$ is not steep then it must be shallow, i.e.\ $a\gs(e-1)l$.  In fact, since $e\mid h_{ij}(\la)=a+l+1$, we find that $a\gs(e-1)l+e-1$.  The definition of $l$ implies that $\la_{i+l+1}<j = \la_i-a$, so
\[\la_i-\la_{i+l+1}\ >\ a\ \gs\ (e-1)(l+1),\]
which implies that for some $k\in\{i,\dots,i+l\}$ we have $\la_k-\la_{k+1}\gs e$.  But this contradicts the assumption that $i>s(\la)$.
\end{pf}

Now we define an operator $S$ on L-partitions.  Suppose $\la$ is an L-partition, and let $s=s(\la)$.  Define
\[S\la = (\la_1-e+1,\la_2-e+1,\dots,\la_s-e+1,\la_{s+2},\la_{s+3},\dots).\]

Note that if $\la$ is an $e$-restricted L-partition, then $S\la=R\la$.  In general, we need to know that $S$ maps L-partitions to L-partitions, in order to allow an inductive proof of Conjecture \ref{slconj}.

\begin{lemma}\label{ssss}
If $\la$ is an L-partition, then so is $S\la$.
\end{lemma}

\begin{pf}
Suppose $\la$ is an L-partition, and that $(i,j)\in[S\la]$.

If $i>s(\la)$, then $(i+1,j)\in[\la]$, and we have
\[a_{ij}(S\la) = a_{(i+1)j}(\la),\qquad l_{ij}(S\la) = l_{(i+1)j}(\la).\]
So if $e\mid h_{ij}(S\la)$, then $e\mid h_{(i+1)j}(\la)$; so by Lemma \ref{shallowsteep}(1) $H_{(i+1)j}(\la)$ is steep, and therefore $H_{ij}(S\la)$ is steep.

Next suppose $i\ls s(\la)$ and $j>\la_{s+1}$.  Then $(i,j+e-1)\in[\la]$ and $a_{ij}(S\la)=a_{i(j+e-1)}(\la)$, $l_{ij}(S\la)=l_{i(j+e-1)}(\la)$.  So if $e\mid h_{ij}(S\la)$, then $e\mid h_{i(j+e-1)}(\la)$, and so $H_{i(j+e-1)}(\la)$ is shallow, and hence $H_{ij}(S\la)$ is shallow.

Finally, suppose that $i\ls s(\la)$ and $j\ls\la_{s+1}$.  Then $(i,j)\in[\la]$, and we have
\[a_{ij}(S\la)=a_{ij}(\la)-e+1,\qquad l_{ij}(S\la)=l_{ij}(\la)-1.\]
So if $e\mid h_{ij}(S\la)$, then $e\mid h_{ij}(\la)$, and hence $H_{ij}(\la)$ is either shallow or steep.  If it is shallow, then we have
\[a_{ij}(S\la)\ =\ a_{ij}(\la)-e+1\ \gs\ (e-1)l_{ij}(\la)-e+1\ =\ (e-1)l_{ij}(S\la),\]
so that $H_{ij}(S\la)$ is shallow.  On the other hand, if $H_{ij}(\la)$ is steep, then
\[l_{ij}(S\la)\ =\ l_{ij}(\la)-1\ \gs\ (e-1)a_{ij}(\la) - 1\ >\ (e-1)a_{ij}(S\la)\]
so $H_{ij}(S\la)$ is steep.
\end{pf}

\begin{eg}
Suppose $e=3$, and let $\la=(9,5,2,1^5)$.  Then we have $s(\la)=2$, so that $S\la=(7,3,1^5)$.  We see that both $\la$ and $S\la$ are L-partitions from the following diagrams.
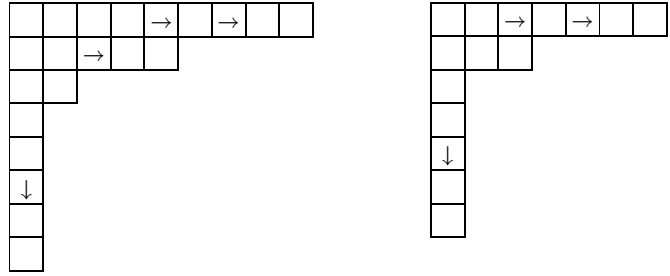
\[\young(\ \ \ \ \shortrightarrow\ \shortrightarrow\ \ ,\ \ \shortrightarrow\ \ ,\ \ ,\ ,\ ,\shortdownarrow,\ ,\ )\qquad\qquad\raisebox{12.7pt}{\young(\ \ \shortrightarrow\ \shortrightarrow\ \ ,\ \ \ ,\ ,\ ,\shortdownarrow,\ ,\ )}\]
\end{eg}

Now we examine the relationship between the operator $S$ and $e$-regularisation.

\begin{lemma}\label{srow}
Suppose $\la$ is an L-partition.  Then
\[GTS\la = CGT\la.\]
\end{lemma}

\begin{pf}
We use induction on $s(\la)$.  In the case $s(\la)=0$ both $\la$ and $S\la=R\la$ are $e$-restricted, i.e.\ $T\la$ and $TS\la$ are $e$-regular, and so $GTS\la = TS\la = TR\la = CT\la = CGT\la$.

Now suppose $s(\la)>0$.  Then $s(R\la)=s(\la)-1$, so we may assume that the result holds with $\la$ replaced by $R\la$.  Put $\zeta=GCT\la$; then by the inductive hypothesis $GTSR\la=CGTR\la=C\zeta$.  Let $\xi$ and $\eta$ be as defined in Lemma \ref{simplereg}, with $x=\la_1$.  Note that
\[x\ =\ \la_1\ \gs\ \la_2+e-1\ =\ l(CT\la)+e-1\ \gs\ l(GCT\la)+e-1\ =\ l(\zeta)+e-1,\]
as required by Lemma \ref{simplereg}.
\clam
$GT\la=G\xi$.
\prof
We have $l(T\la)=\la_1=l(\xi)$ and $GCT\la=\zeta=C\xi$, and Lemma \ref{cggc} gives the result.
\malc
\clam
$GTS\la=G\eta$.
\prof
Since $s(\la)>0$, $S\la$ may be obtained from $SR\la$ by adding a row of length $\la_1-e+1$; hence $TS\la$ may be obtained from $TSR\la$ by adding a column of length $\la_1-e+1$.  So we have $l(TS\la)=\la_1-e+1=l(\eta)$, and
\[GCTS\la = GTSR\la = C\zeta = C\eta,\]
and again we may appeal to Lemma \ref{cggc}.
\malc
Now Lemma \ref{simplereg} combined with these two claims gives the result.
\end{pf}

Next we prove a simple lemma which gives an equivalent statement to the condition $MG\la=GT\la$ in the presence of a suitable inductive hypothesis.

\begin{lemma}\label{equiv}
Suppose $\la$ is an L-partition, and that $MG\mu=GT\mu$ for all L-partitions $\mu$ with $|\mu|<|\la|$.  Then $MG\la=GT\la$ if and only if $GS\la=JG\la$.
\end{lemma}

\begin{pfnb}
Since $|G\la|=|GT\la|$, we have
\begin{align*}
MG\la = GT\la \quad&\Longleftrightarrow\quad MJG\la = CGT\la\tag*{by Proposition \ref{xuv}}\\
\quad&\Longleftrightarrow\quad MJG\la = GTS\la\tag*{by Lemma \ref{srow}}\\
\quad&\Longleftrightarrow\quad MJG\la = MGS\la\tag*{by the inductive hypothesis and Lemma \ref{ssss}}\\
\quad&\Longleftrightarrow\quad JG\la=GS\la.\tag*{\qedsymbol}
\end{align*}
\end{pfnb}

We now require one more lemma concerning the regularisations of L-partitions.

\begin{lemma}\label{sssreg}
Suppose $\la$ is an L-partition with $s(\la)>0$ and $\la_1\gs l(\la)$.  Then:
\begin{enumerate}
\item
$(G\la)_1=\la_1$;
\item
$(G\la)_1-(G\la)_2\gs e-1$;
\item
$(GS\la)_1=(S\la)_1$.
\end{enumerate}
\end{lemma}

\begin{pfenum}
\item
Obviously $(G\la)_1\gs\la_1$, so it suffices to show that $[\la]$ does not contain a node in ladder $(e-1)\la_1+1$. If it does, let $(i,j)$ be the rightmost such node.  Since $(i,j)\neq(1,\la_1+1)$, we have $i\gs e$ and we know that the node $(i-e+1,j+1)$ does not lie in $\la$; in other words, $(T\la)_j-(T\la)_{j+1}\gs e$.  This means that $j\ls t(\la)$, and so by Lemma \ref{gse-1} we have $i\ls l(\la)-(e-1)(j-1)$, so that
\[l(\la)\gs i+(e-1)(j-1)=(e-1)\la_1+1>\la_1,\]
contrary to hypothesis.
\item
By part (1), we must show that $(G\la)_2\ls \la_1-e+1$, i.e.\ that $[\la]$ does not contain a node in ladder $2+(e-1)(\la_1-e+1)$.  Supposing otherwise, we let $(i,j)$ be the rightmost such node.  Arguing as above, we find that
\[\la_1\gs l(\la)\gs i+(e-1)(j-1) = 2+(e-1)(\la_1-e+1),\]
and this rearranges to yield $\la_1<e$, which is absurd given that $s(\la)>0$.
\item
Obviously $(GS\la)_1\gs(S\la)_1=\la_1-e+1$, so it suffices to show that $[S\la]$ does not contain a node in ladder $1+(e-1)(\la_1-e+1)$.  Arguing as above, such a node would have to be of the form $(i,j)$ with $j\ls t(S\la)\ls t(\la)$.  But then $(TS\la)_j=(T\la)_j-1$, so $[\la]$ contains the node $(i+1,j)$, which lies in ladder $2+(e-1)(\la_1-e+1)$.  But it was shown in (2) that this is not possible.
\end{pfenum}

\begin{pfof}{Conjecture \ref{slconj} (`if' part)}
We proceed by induction on $|\la|$.  It is clear that $\la$ is an L-partition if and only if $T\la$ is, so Conjecture \ref{slconj} holds for $\la$ if and only if it holds for $T\la$.  If either $\la$ or $T\la$ is $e$-regular, then the result follows from Theorem \ref{boxthm}, so we assume that $\la$ is neither $e$-regular nor $e$-restricted; in particular, $s(\la)>0$.  By replacing $\la$ with $T\la$ if necessary, we assume also that $\la_1\gs l(\la)$.

%
\clam
$(JG\la)_1=\la_1-e+1$, and $RJG\la=JGR\la$.
\prof
This follows from Lemma \ref{sssreg}(1--2), given the definition of the operator $J$.
\malc

\clam
$(GS\la)_1=\la_1-e+1$, and $RGS\la=GRS\la$.
\prof
We have $(S\la)_1=\la_1-e+1$ by definition, and $(GS\la)_1=(S\la)_1$ by Lemma \ref{sssreg}(3).  The second statement follows from Lemma \ref{rggr}.
\malc

By induction (replacing $\la$ with $R\la$) we have $MGR\la=GTR\la$, and by Lemma \ref{equiv} (and the inductive hypothesis) this gives $JGR\la=GSR\la$.  Since obviously $GSR\la=GRS\la$, the two claims yield $JG\la=GS\la$.  Now applying Lemma \ref{equiv} again gives the result.
\end{pfof}

\section{The Fock space and $v$-decomposition numbers}\label{onlyif}

In this section, we complete the proof of Conjecture \ref{slconj} using $v$-decomposition numbers.  We give only a very brief sketch of the background material needed, since this is discussed at length elsewhere; in particular, the article of Lascoux, Leclerc and Thibon \cite{llt} is an invaluable source.

Fix $e\gs2$, let $v$ be an indeterminate over $\bbq$, and let $\calu$ be the quantum algebra $U_v(\widehat{\mathfrak{sl}}_e)$ over $\bbq(v)$.  There is a module $\calf$ for this algebra called the \emph{Fock space}, which has a  \emph{standard basis} indexed by (and often identified with) the set of all partitions.  The submodule generated by the empty partition is isomorphic to the \emph{basic representation} of $\calu$.  This submodule has a \emph{canonical} $\bbq(v)$-basis
\[\big\{G(\mu)\ \big|\ \mu\text{ an $e$-regular partition}\big\}.\]
The \emph{$v$-decomposition numbers} are the coefficients obtained when the elements of the canonical basis are expanded in terms of the standard basis, i.e.\ the coefficients $d_{\la\mu}(v)$ in the expression
\[G(\mu) = \sum_{\la}d_{\la\mu}(v)\la.\]

We shall need to quote two results concerning $v$-decomposition numbers; one concerning the Mullineux map, and the other concerning $e$-regularisation.  The first of these involves the $e$-weight of a partition, defined in \S\ref{hooks}.

\begin{thm}\label{lltthm}\thmcite{llt}{Theorem 7.2}
Suppose $\la$ and $\mu$ are partitions with $e$-weight $w$, and that $\mu$ is $e$-regular.  Then
\[d_{(T\la)(M\mu)}(v) = v^wd_{\la\mu}(v^{-1}).\]
\end{thm}

The second result we need requires a definition.  Given a partition $\la$, let $z(\la)$ be the number of nodes $(i,j)\in[\la]$ such that $e\mid h_{ij}(\la)$ and $H_{ij}(\la)$ is steep.  Now we have the following result.

\begin{thm}\label{regthm}\thmcite{mfreg}{Theorem 2.2}
For any partition $\la$,
\[d_{\la(G\la)}(v) = v^{z(\la)}.\]
\end{thm}

\begin{rmk}
Note that in \cite{mfreg} an alternative convention for the Fock space is used: our $d_{\la\mu}(v)$ is written in \cite{mfreg} as $d_{(T\la)(T\mu)}(v)$.  Accordingly, the statement of \cite[Theorem 2.2]{mfreg} involves shallow hooks rather than steep hooks.  We hope that no confusion will result.
\end{rmk}

Now we combine these theorems.  First we note the following obvious result about $e$-weight and the function $z$.

\begin{lemma}\label{zlemma}
Suppose $\la$ is a partition with $e$-weight $w$.  Then $T\la$ also has $e$-weight $w$, and $z(T\la)$ equals the number of nodes $(i,j)\in[\la]$ such that $e\mid h_{ij}(\la)$ and $H_{ij}(\la)$ is shallow.  Hence $\la$ is an L-partition if and only if $w=z(\la)+z(T\la)$.
\end{lemma}

Now we can complete the proof of Conjecture \ref{slconj}.

\begin{pfof}{Conjecture \ref{slconj} (`only if' part)}
Suppose $MG\la=GT\la$, and that $\la$ has $e$-weight $w$.  Then we have
\begin{align*}
v^{z(T\la)} &= d_{(T\la)(GT\la)}(v)\tag*{by Theorem \ref{regthm}}\\
&= d_{(T\la)(MG\la)}(v)\tag*{by hypothesis}\\
&= v^wd_{\la(G\la)}(v^{-1})\tag*{by Theorem \ref{lltthm}}\\
&= v^w.v^{-z(\la)}\tag*{by Theorem \ref{regthm}}
\end{align*}
so that $w=z(\la)+z(T\la)$.  Now Lemma \ref{zlemma} gives the result.
\end{pfof}

\end{document}